\documentclass{article}
\usepackage{cite}

 \newcommand {\theoremstyle} [1] { }

 \theoremstyle{plain}
 \theoremstyle{definition}

 \theoremstyle{remark}

\usepackage[pagewise]{lineno}%\linenumbers
\usepackage {amssymb}

\usepackage {amssymb}
\usepackage{amsmath,amsfonts, color}
\usepackage{graphicx}

\makeatother

\author{Colin Rogers$^1$ and Pablo Amster$^2$}

\title{A  2+1-dimensional extended Dym equation: moving boundary problems solvable via Painlev\'e II symmetry reduction%
}
\date{}

\begin{document}
\maketitle

\begin{center}
$^1$ {School of Mathematics and Statistics}\\
   {The University of New South Wales}\\
    {Sydney, NSW 2052, Australia}
    \end{center}

\begin{center}
$^2$ Departamento de Matem\'atica, \\
Facultad de Ciencias Exactas y Naturales\\
Universidad de Buenos Aires
and 
IMAS - CONICET\\
Ciudad Universitaria, Pabell\'on I,
(1428) Buenos Aires, Argentina 
\end{center}

\begin{center}
 c.rogers@unsw.edu.au --- pamster@dm.uba.ar
\end{center}
\bigskip

%\begin{document}

%\runningtitle{A novel 2+1-dimensional extended Dym equation}        % if too long for running head

% \and
%				Luigi Vergori



\maketitle

\begin{abstract}
A  2+1-dimensional extension of the solitonic Dym equation is shown to admit a Painlev\'e II symmetry reduction which permits the exact solution of a class of Stefan-type moving boundary problems.
 \medskip
 
 \noindent \textbf{Keywords}: Extended Dym Equation; Ermakov- Painlev\'e II; Moving Boundary Problem.

 \smallskip
 \noindent \textbf{MSC2020}: {80A22, 35Q51}
\end{abstract}

\section{Introduction}
In 1+1-dimensional soliton theory, it was established in \cite{cr15} that novel classes of nonlinear moving boundary problems of Stefan-type and their reciprocal associates for the canonical Dym equation are amenable to exact solution via Painlev\'e II symmetry reduction. In a subsequent development \cite{cr17}, moving boundary problems for an integrable extension of the Dym equation with application in peakon hydrodynamics (Camassa-Holm \cite{rcdh93}) were solved by linkage to an associated class of Stefan-type problems for the standard Dym equation. Exact solution was thereby derived notably in terms of Yablonski-Vorob'ev polynomials via Painlev\'e II reduction. It is remarked that in \cite{wscr99} both the canonical 1+1-dimensional Dym equation and its Camassa-Holm type extension have been embedded in a class of solitonic torsion evolution equations associated with binormal motions of an inextensible curve and which is reciprocal related (\cite{jkcr82}) to the solitonic m$^2$KdV equation originally derived in \cite{af80} and subsequently in an independent study by Calogero and Degasperis \cite{fcad85}.

A 2+1-dimensional integrable version of the Dym equation, namely
\begin{equation} \label{1}
r_t+r^3r_{xxx}+(3/r)(r^2\partial^{-1}_x\left(\frac{r_y}{r^2}\right)_y)=0 \end{equation}
was originally introduced by Konopelchenko and Dubrovsky in \cite{bkvd84} and subsequently linked in \cite{cr87} via a 2+1-dimensional reciprocal-type transformation \cite{cr86} to a singularity manifold equation which results from Painlev\'e integrability criteria applied to the canonical Kadomtsev-Petviashvili equation of shallow water hydrodynamics \cite{bkvp70}. In a recent development \cite{bkcrpa26} an integrable extension of \eqref{1} has been derived i.e.
\begin{equation} \label{2}
u_t+2\partial_x(1-\partial_{xx})\left(\frac{1}{u^{1/2}}\right)+6u^2[\ u^{-1}\partial^{-1}_x(u^{1/2})_y\ ]_y=0 \end{equation}
which is amenable to a version of the $\bar{\partial}$-dressing method of Zakharov and Shabat \cite{vzas74} with extensions detailed in \cite{bk92}. In \cite{gwztxgjl19} another kind of %novel
2+1-dimensional extension of the canonical solitonic Dym equation has recently been set down which is equivalent to
\begin{equation} \label{3}
\partial p/\partial t=2[\ (p^{-1/2})_{xxx}+(p^{-1/2})_{yyy}+(p^{-1/2})_{xyy}+(p^{-1/2})_{xxy}\ ].\end{equation}
The latter is shown here to admit a Painlev\'e II symmetry reduction which allows exact solution of a class of nonlinear moving boundary problems of Stefan-type.
It is remarked that a 2+1- dimensional deformed Dym equation cognate to (3) had previously been introduced in \cite{feran}.

\section{Painlev\'e II Symmetry Reduction}

Here, a symmetry reduction of the extended 2+1-dimensional Dym equation is sought via the ansatz
\begin{equation} \label{4}
p^{-1/2}=(t+a)^mP\left(\frac{x+\alpha^*y}{(t+a)^n}\right), \quad \alpha^* \in  \mathbb{R} \end{equation}
whence insertion of the latter into \eqref{3} yields
\begin{equation} \label{5}
\begin{array}{l} -2m(t+a)^{-2m-1}P^{-2}-n\xi(t+a)^{-2m-1}(-2)P^{-3}P' \\[2mm]
\quad =2[\ (t+a)^m(t+a)^{-3n}P^{'''}+(t+a)^m(t+a)^{-3n}\alpha^{*3}P^{'''} \\[2mm]
\qquad +(t+a)^m(t+a)^{-3n}\alpha^{*2}P^{'''}+(t+a)^m(t+a)^{-3n}\alpha^*P^{'''}]\end{array}\end{equation}
with $\xi=(x+\alpha^* y)/(t+a)^n$.

Accordingly,
\begin{equation} \label{6}
P^{'''}[\ 1+\alpha^*+\alpha^{*2}+\alpha^{*3}\ ]=-mP^{-2}+n\xi P^{-3}P', \end{equation}
together with
\begin{equation} \label{7}
m=(3n-1)/3 \end{equation}
and integration of \eqref{6} shows that
\begin{equation} \label{8}
\begin{array}{l} [\ PP^{''}-P^{'2}/2\ ]\ [\ 1+\alpha^*+\alpha^{*2}+\alpha^{*3}\ ] \\[2mm]
\quad = -\dfrac{n\xi}{P}-(m-n)\displaystyle\int \dfrac{1}{P} d\xi+\mathrm{I}\ , \quad \mathrm{I} \in  \mathbb{R}. \end{array}\end{equation}
On introduction of
\begin{equation} \label{9}
w=aP_\xi \end{equation}
together with
\begin{equation} \label{10}
s=[\ PP^{''}-P^{'2}/2\ ]\ [\ 1+\alpha^{*}+\alpha^{*2}+\alpha^{*3}\ ]+n\xi/P=\mathrm{I}-(m-n)\displaystyle\int\frac{1}{P} d\xi \end{equation}
there results
\begin{equation} \label{11}
w_s=-\frac{aP^{''}P}{m-n}, \end{equation}
\begin{equation} \label{12}
w_{ss}=d/\partial\xi\left[-\frac{aP^{''}P}{m-n}\right]\frac{d\xi}{ds}=\frac{a}{(m-n)^2}[\ P^{'''}P+P'P^{''}\ ]P. \end{equation}
On use of the relations \eqref{6} and \eqref{10} respectively
\begin{equation} \label{13}
P^{'''}=[\ -mP^{-2}+n\xi P^{-3}P'\ ]/[\ 1+\alpha^*+\alpha^{*2}+\alpha^{*3}\ ], \end{equation}
\begin{equation} \label{14}
PP^{''}=[\ s-n\xi P^{-1}\ ]/[\ 1+\alpha^*+\alpha^{*2}+\alpha^{*3}\ ]+P^{'2}/2. \end{equation}
with insertion into \eqref{12} there results
\begin{equation} \label{15}
w_{ss}=\dfrac{a}{(m-n)^2}\left[\dfrac{-m+n\xi\ P^{-1}P'+P'(s-n\xi P^{-1})}{1+\alpha^*+\alpha^{*2}+\alpha^{*3}}+\frac{P^{'3}}2\right]. \end{equation}

The latter with $s=\epsilon z$ and the appropriate scaling results in the canonical Painlev\'e II equation
\begin{equation} \label{16}
w_{zz}+2w^3+zw+\alpha=0, \end{equation}
with parameter
\begin{equation} \label{17}
\alpha=-ma/\left[(m-n)^2(1+\alpha^*+\alpha^{*2}+\alpha^{*3})\right]. \end{equation}

\section{A Class of Moving Boundary Problems}

Here, moving boundary problems for the extended 2+1-dimensional Dym equation \eqref{6} in a region $0<x+\alpha^*y<S(t)$ are shown to admit exact solution via the Painlev\'e II symmetry reduction with $w=z^{-1}$, $\alpha=-1$. The boundary conditions for the Stefan-type moving boundary problem under consideration are given by
\begin{equation} \label{18}
\begin{array}{c} -2[\ (p^{-1/2})_{xx}+(p^{-1/2})_{xy}+(p^{-1/2})_{yy}+\partial^{-1}_x(p^{-1/2})_{yyy}\ ]=L_mS^i\dot{S}\ , \\[3mm]
t>0 \end{array}\end{equation}
\begin{equation} \label{19}
p=P_mS^j\ , \quad t>0 \end{equation}
on $x+\alpha^*y=S(t)$ together with
\begin{equation} \label{20}
\begin{array}{l} -2[\ (p^{-1/2})_{xx}+(p^{-1/2})_{xy}+(p^{-1/2})_{yy}+\partial^{-1}_x(p^{-1/2})_{yyy}\ ]|_{x+\alpha^*y=0} \\[4mm]
\qquad\qquad\qquad =H_0(t+a)^k\ , \quad t>0 \end{array}\end{equation}
wherein, in the preceding $S(t)=\gamma(t+a)^n$. $P(\xi)$ in the symmetry reduction representation \eqref{4} is here determined by the triad of relations
\begin{equation} \label{21}
\begin{array}{c} w=aP'(\xi)\ , \quad w=-1/z \\[2mm]
\epsilon z=\mathrm{I}-(m-n)\displaystyle\int\dfrac{1}{P}d\xi \end{array}\end{equation}
whence
\begin{equation} \label{22}
P^{''}/P^{'}=-\frac{a(m-n)}{\epsilon}P'/P. \end{equation}
On integration, the latter yields
\begin{equation} \label{23}
P'P^{a(m-n)/\epsilon}=\sigma_1\ \in \mathbb{R} \end{equation}
whence
\begin{equation} \label{24}
P=[\ \sigma_1 [\ a(m-n)/\epsilon+1\ ]\xi+\sigma_2\ ]^{1/[a(m-n)/\epsilon+1]}, \quad \sigma_2\ \in \mathbb{R}. \end{equation}

$$\textbf{Boundary Conditions}$$
I. \ $-2[\ (p^{-1/2})_{xx}+(p^{-1/2})_{xy}+(p^{-1/2})_{yy}+\partial^{-1}_x(p^{-1/2})_{yyy}\ ]=L_m S^i\dot{S}$ \\[-2mm] $$\text{on} \quad x+\alpha^*y=S(t)=\gamma(t+a)^n$$

whence $\xi=\gamma$. \\

Insertion of the Painlev\'e II symmetry reduction requirement yields
\begin{equation} \label{25}
-2(t+a)^{m-2n}[\ 1+\alpha^*+\alpha^{*2}+\alpha^{*3}\ ] P^{''}(\gamma)=n\gamma^{i+1}(t+a)^{ni+n-1}L_m \end{equation}
so that $i=-2m/n$ by virtue of the relation \eqref{7} while
\begin{equation} \label{26}
L_m=-2[\ 1+\alpha^*+\alpha^{*2}+\alpha^{*3}\ ]P^{''}(\gamma)/n\gamma^{i+1}.\end{equation}
II. \ $p=P_mS^j$ \quad on \quad $x+\alpha^*y=S(t), \quad t>0$.\\

This condition yields
\begin{equation} \label{27}
P_m=\gamma^{-j}P^{-2}(\gamma) \end{equation}
together with $j=-2m/n$.\\[3mm]
III. \ $-2[\ (p^{-1/2})_{xx}+(p^{-1/2})_{xy}+(p^{-1/2})_{yy}+\partial^{-1}_x(p^{-1/2})_{yyy}\ ]|_{x+\alpha^*y=0}$ \\[-2mm] $$\qquad\qquad =H_0(t+a)^k, \quad t>0.$$

Here, the Painlev\'e II symmetry reduction requires that $k=m-n$ together with
\begin{equation} \label{28}
H_0=-2[\ 1+\alpha^*+\alpha^{*2}+\alpha^{*3}\ ]P^{''}(0). \end{equation}

In conclusion, it is remarked that the present analysis may be extended to generalisations of the 2+1-dimensional equation (3) with either temporal or spatial modulation
as determined by the application of a class of involutary transformations with genesis in classical Ermakov theory \cite{16} .The procedure in the case of temporal modulation is described in [12] and in the spatially modulated application to coupled solitonic systems of sine-Gordon, de Moulin  and Manakov-type in \cite{17}.

%\begin{acknowledgements}
%\end{acknowledgements}



\end{document}